\documentclass[a4paper,10pt]{article}
\usepackage[cp1251]{inputenc}
\usepackage[english]{babel}
\usepackage[T2A]{fontenc}
\usepackage{mathtext}
\usepackage{amsmath,amssymb,amsthm,amsfonts}
\usepackage{indentfirst}

\DeclareSymbolFont{T2Aletters}{T2A}{cmr}{m}{it}

\newtheorem{defin}{Definition}
\newtheorem{tver}{Proposition}
\newtheorem{nasl}{Corollary}
\newtheorem{umova}{Assumption}
\newtheorem{remark}{Remark}
\newtheorem{teorema}{Theorem}
\newtheorem{lema}{Lemma}
\begin{document}

\textit{Transactions of NAS of Azerbaijan}, (2010), vol. XXX, No 4, pp. 141-152.

\vspace{3mm}
\begin{footnotesize}
Received February 02, 2010; Revised May 11, 2010.
\end{footnotesize}
\begin{center}
\Large{\textbf{KILLED MARKOV DECISION PROCESSES ON FINITE TIME INTERVAL FOR COUNTABLE MODELS}}\footnote{
Revised and corrected version of the paper published in \textit{Transactions of NAS of Azerbaijan}, (2010), vol. XXX, No 4, pp. 141-152.
}
\end{center}

\vspace{3mm}

\begin{center}
\textbf{Nestor R. PAROLYA and Yaroslav I. YELEYKO}
\end{center}

\vspace{3mm}

\begin{abstract}
\textit{We consider killed Markov decision processes for countable models on a finite time-interval.
Existence of a uniform $\varepsilon$-optimal policy is proven. We show the correctness of the fundamental equation.
The optimal control problem is reduced to a similar problem for the derived model. We receive an optimality
equation and a method for the construction of simple optimal policies. The sufficiency of simple policies
for countable models is proven. We show the correctness of the Markovian property. Additionally, a dynamic
programming principle is considered.}
\end{abstract}

\vspace{3mm}

\textit{Classification:~} 90C40.


\paragraph{1. Introduction.}
Markov decision processes arise in the different areas of the economics,
in particular for the economic work planning of the separate business, economic sector or 
entire economics. At the beginning of each period we can build a plan for the next period  
knowing the last achieved state. The system development can be described mathematically
as a deterministic process if we assume that the position of the system at the end of each period is
uniquely defined by the state at the end of the period and by a plan for this period.

It is necessary to consider the influence of such factors as meteorological conditions, 
demographic transition, demand fluctuations, the imperfection of the compound 
production processes coordination, scientific discoveries and inventions etc. Stochastic models take into account these factors: if we know the state at the beginning of the period and the plan, we can only calculate the probability distribution for the next period. Therefore, 
leaving aside the system states in the past periods we come to the idea of Markov decision process
("the future depends not on the past, but only on the present").

The Markov decision processes are well described in [1]:
the definition of Markov decision process is given, the concept of "model" $Z^\mu$ is presented,
the definition of policy $\pi$ is given, the assessment of policy - $\omega(\pi)$ and $\nu$ - assessment
of process $Z^\mu$ are defined, the existence of a uniform $\varepsilon$-optimal policy is proved,
the optimality equation and method for simple optimal policies constructing are presented, the sufficient of simple 
policies for countable models is proved, the correctness of the Markovian property is shown and dynamic programming 
principle is considered.

In [1] the model does not take into account one risk factor, namely the probability of
bankruptcy at some determined moment of time. As a result, we come to the idea of
killed Markov decision process where the business can crash with some nonzero probability at every moment of time,
with the exception of the initial state.

The concept of the killed Markov decision process brings us closer to the real economic system
which is not common without risk.

\paragraph{2. Killed Markov decision process.}

Let $X_t (t=m,\ldots,n)$ and let $A_t (t=m+1,\ldots,n)$ be countable or finite sets and at least one of them is countable. 
To the arbitrary $a\in A_t$ is assigned a probability distribution
$p(\cdot|a)=\mathbb{P}(x_t=x|a_t=a,x_{t-1})$ on $X_t$.

\begin{defin}
The function $p$ which defines the law of the transition from $A_t$ to $X_t$ is called the
 \textbf{transition function}.
\end{defin}

\begin{defin}
The point $x^*=x_m \in X_t$ is called \textbf{killed state}, and
$p(x^*|a)$ - the \textbf{probability of kill} if
$\mathbb{P}(x_{t+1}=x^*|a_t=a)=\mathbb{P}(x_{t+1}=x_m|a_t=a)\equiv
p(x^*|a), x_m \in X_m$.
\end{defin}

\begin{remark} In other words, the system moves into the initial(home) state when it hits a killed state(process is killed).
\end{remark}
From the definition of the killed state it follows:
$$\forall a \in A_t \ \exists x^* \in X_t: p(x^*|a)=1-\sum_{x \in X_t\setminus
x^*}{p(x|a)>0}.$$

\begin{defin}[Killed Markov decision process]
A killed Markov decision process on a time interval $[m,n]$ is defined through the following objects:

1. Sets $X_m,\ldots,X_n$(spaces of states);

2. Sets $A_{m+1},\ldots,A_n$(spaces of actions);

3. The projection mapping $j:A\rightarrow X$ where
$A=\bigcup\limits_{t=m+1}^{n}A_t$, $X=\bigcup\limits_{t=m}^{n}X_t$:
$j(A_t)=X_{t-1}\setminus \{x^*\},x^* \ \in X_{t-1},
(t=m+2,\ldots,n)$ and $j(A_{m+1})=X_m$;

4. The probability distribution
$p(\cdot|a)=\mathbb{P}(x_{t}=x|a_t=a,x_{t-1})$ on $X_t$ with killed states
$$\mathbb{P}(x_{t+1}=x^*|a_t=a)=\mathbb{P}(x_{t+1}=x_m|a_t=a)\equiv
p(x^*|a)>0;$$

5. The function $q$ on $A$ (reward function);

6. The function $r$ on $X_n$  (terminal reward);

7. The function $c$ (crash function), defined on the killed states
$c(x^*)=-\sum\limits_{i=m+1}^t\max\limits_{a_i \in A_i}q(a_i), x^*
\in X_t, t=m+1,\ldots,n$ (function $c$ ensures a total bankruptcy - total loss of accumulated capital or more);

8. The initial distribution $\mu$ on $X_m$.

A stochastic process defined through (1-8) is called the \textbf{\textit{killed Markov decision process}} or the \textbf{\textit{model}} and it is
denoted by $Z_{\mu}^*$. If the initial distribution $\mu$ is concentrated at the point $x$, we shall write $Z_{x}^*$.
\end{defin}

\begin{defin}
The trajectory $l=x_{m}a_{m+1}x_{m+1}\ldots a_{n}x_{n}$ is called the \textbf{way}.
The set of all ways we denote $L=X\times(X\times A)^n$.
\end{defin}

Our goal is to find a decision method which maximizes the mathematical expectation of the assessment of way $l$:
$$I(l,x^*)=\sum\limits_{t=m+1}^{n}[q(a_t)+c(x_t^*)]+r(x_n),\eqno
(2.1)$$ where:

$x^*=(x_{m+1}^*,\ldots,x_n^*)$ - vector of killed states;

$l=x_{m}a_{m+1},\ldots,a_{n}x_{n}$ - way.

\vspace{4mm}

The decision method is meant to be some \textit{policy}.

\paragraph{3. Policies.} 

\begin{defin}
Let $A(x)\subset A$ is the set of all available actions at the state $x\in X$. 
$\varphi(x):X\rightarrow A(x)$ is called the \textbf{simple
policy} if $\varphi(x_{t-1})=a_t$ for arbitrary $x_t$ which is not a killed state with the probability distribution
$p(\cdot|a_t)(m<t\leq n)$ and $x_m$ with the initial distribution $\mu$.
\end{defin}

\begin{remark} When we use the simple policy $\varphi(x)$ we get the way $l=x_{m}a_{m+1},\ldots,a_{n}x_{n}$.
\end{remark}
\begin{defin}
The mapping $\pi:H\rightarrow \pi(\cdot|h \in H)$ is called a \textbf{killed policy}, where $\pi(\cdot|h
\in H)$ is a probability distribution on $A(x_{t-1})$ and $H=X\times(A\times X)^{t-1}$ is a space of histories up to epoch $m\leq t-1\leq n$ $(h \in H\Leftrightarrow h=x_{m}a_{m+1},\ldots,a_{t-1}x_{t-1})$.
\end{defin}

\begin{remark} 
Obviously, $x_{t-1}\neq x^*$.
\end{remark}
\begin{defin}
Killed policy $\pi(\cdot|h)$ is called a
\textbf{Markov policy} if $\pi(\cdot|h)=\pi(\cdot|x_{t-1})$.
\end{defin}

The next conceptions can not be well-defined without the assumption:
\begin{umova}
The reward function $q$ and the terminal reward function $r$ have the \textbf{supremum,}
$\exists \sup\limits_{a\in A} q(a)$ and $\exists \sup\limits_{x\in X_n} r(x)$.
\end{umova}
\begin{defin}
Let $p(\cdot|a)$ be the transition function and let
$\pi(\cdot|h)$ be a policy. Every initial distribution $\mu$ is assigned to a probability distribution $P^*$ in the space $L$ which has such the notation:
$$P^*(l,x^*)=P^*(x_{m}a_{m+1},\ldots,a_{n}x_{n},~x_{m+1}^*,\ldots,x_{n}^*)=$$
$$=\mu(x_m)\pi(a_{m+1}|x_m)p(x_{m+1}|a_{m+1})p(x_{m+1}^*|a_{m+1})
\cdot\ldots\cdot\pi(a_n|h_{n-1})p(x_n|a_n)p(x_{n}^*|a_n) \eqno
(3.1)$$
\end{defin}

\begin{remark} 
After the definition of the measure $P^*$ the way $l$ can be interpreted as a stochastic process. Additionally, this process is called the Markov process if the policy $\pi$ is a Markov policy.
\end{remark}
For all functions $\xi$ from space $L$ the mathematical expectation of $\xi$ is given by
$$E^*(\xi)=\sum\limits_{l \in L}\xi(l)P^*(l,x^*) \eqno (3.2)$$

The assessment (2.1) of the way $l$ is an example of such function. Next, we denote its expectation $\omega$:
$$\omega=E^*I(l,x^*)=E^*[\sum\limits_{t=m+1}^{n}[q(a_t)+c(x_t^*)]+r(x_n)] \eqno (3.3)$$

\begin{defin}[Assessment of policy] 
The value $\omega$ from (3.3) is called the \textbf{assessment of policy} $\pi$ and is the function of the variable $\pi$ ($\omega=\omega(\pi)$) for the killed Markov decision process $Z_{\mu}^*$.

\end{defin}

The goal of the research is the maximization of function $\omega(\pi)$.

\begin{defin}[Assessment of process]
$\nu\equiv\sup\limits_{\pi}\omega(\pi)$ is called the \textbf{assessment of killed Markov decision process} $Z_{\mu}^*$ or \textbf{assessment of initial distribution} $\mu$.
\end{defin}
\begin{remark} $\nu(x^*)=c(x^*)$.
\end{remark}
\begin{defin}[$\varepsilon$-optimal policy]
A killed policy $\pi$ is called \textbf{$\varepsilon$-optimal} for $Z^*_\mu$ if $\forall \varepsilon>0:\
\omega(\mu,\pi)\geq\nu(\mu)-\varepsilon$.
\end{defin}
\begin{defin}[Uniform $\varepsilon$-optimal policy]
A killed policy is called \textbf{uniform} $\varepsilon$-\textbf{optimal} or $\varepsilon$-\textbf{optimal for process}
$Z^*$ if $\pi$ is $\varepsilon$-optimal for $Z_{\mu}^*$ for all
$\mu$ - initial distribution.
\end{defin}

\paragraph{4. Existence of uniform $\varepsilon$-optimal policy.}

Let $\pi_x$ is $\varepsilon$-optimal policy for process $Z^*_x$. Its existence follows from the definition of the supremum.

We want to build a killed policy $\pi$ which is $\varepsilon$-optimal for the model $Z^*$ by using a sequence of the killed policies $\pi_x$.

It's natural to use the policy $\pi_x$ when $x$ is a starting point. Formally,

$$\bar{\pi}(\cdot|h)=\pi_{x(h)}(\cdot|h) \eqno (4.1)$$

where $x(h)$ - the initial state of history $h$. 
It is clear that formula (4.1) defines some policy $\bar{\pi}$ and this policy will be $\varepsilon$-optimal. It means that
$\forall \varepsilon\geq0:\omega(x,\bar{\pi})=\omega(x,\pi_x)\geq\nu(x)-\varepsilon, \forall x
\in X_m$.

\begin{tver}[Existence of the uniform $\varepsilon$-optimal killed policy]
Every killed policy $\bar{\pi}$ from (4.1) which is $\varepsilon$-optimal, i.e.
$$\omega(x,\bar{\pi})\geq\nu(x)-\varepsilon,(x \in X_m),~\forall \varepsilon\geq0$$
is uniform $\varepsilon$-optimal. It means that $\forall \mu,\forall
\varepsilon\geq0: \
\sup\limits_{\pi}\omega(\mu,\pi)\leq\omega(\mu,\bar{\pi})+\varepsilon.$
\end{tver}

\begin{flushleft}
\textbf{Proof}. From (3.1)-(3.3) it follows that $\forall \pi$:
\end{flushleft}

$$\omega(\mu,\pi)=\sum\limits_{l \in L}I(l,x^*)P^*(l,x^*)=\sum\limits_{X_m}\mu(x)\omega(x,\pi). \eqno (4.2)$$

Hence, it appears
$$\omega(\mu,\pi)=\sum\limits_{X_m}\mu(x)\omega(x,\pi)\leq\sum\limits_{X_m}\mu(x)\nu(x)\leq\sum\limits_{X_m}\mu(x)[
\omega(x,\bar{\pi})+\varepsilon]=\omega(\mu,\bar{\pi})+\varepsilon.$$

From the received inequalities it follows that

$$\sup\limits_{\pi}\omega(\mu,\pi)\leq\sum\limits_{X_m}\mu(x)\nu(x),\eqno(4.3)$$
$$\omega(\mu,\bar{\pi})\geq\sum\limits_{X_m}\mu(x)\nu(x)-\varepsilon. \eqno(4.4)$$

According to the arbitrariness of $\varepsilon>0$ we get now from (4.3) and (4.4)
$$\sup\limits_{\pi}\omega(\mu,\pi)=\sum\limits_{X_m}\mu(x)\nu(x)\leq\omega(\mu,\bar{\pi})+\varepsilon. \eqno (4.5)$$

So the policy $\bar{\pi}$ is uniform $\varepsilon$-optimal. \textit{\textbf{Proposition 1 is proved.}}
\begin{flushleft}
\textbf{Corollary 1.}
For all initial distributions $\mu$:
\end{flushleft}

$$\nu(\mu)=\mu\nu. \eqno (4.6)$$
\begin{flushleft}
\textbf{Proof.} It follows from $\nu(\mu)=\sum\limits_{X_m}\mu(x)\nu(x)=\mu\nu.$
\end{flushleft}

\vspace{4mm}

\begin{remark} Formulas (4.2) and (4.6) allow us to reduce the analysis of the processes $Z_\mu^*$ for all $\mu$ to the analysis of the processes $Z^*_x,~\forall
x\in X_m$.
\end{remark}
\vspace{4mm}

The policy $\pi$ is built of the sequence $\pi_x, (x\in X_m)$ and has the following property (1):
\begin{center}
\textit{For all initial distribution of the state $x\in X_m$ the probability distributions in space $L$ which are assigned to the policies $\pi$ and $\pi_x$  from (3.1) are equal.}
\end{center}
\begin{defin}
If $\bar{\pi}$ satisfies the property (1) then $\bar{\pi}$ is called the \textbf{combination of policies} $\pi_x$.
\end{defin}

\paragraph{5. Derived model and fundamental equation.}

The decision process is a quite number of consecutive steps. The first step is the choice of probability distribution on $A_{m+1}$ which depends on initial state. Since the choice is taken every initial distribution $\mu$ on $X_m$ accords with probability distribution $\acute{\mu}$ on $X_{m+1}$.
Now we consider $\acute{\mu}$ as initial distribution in moment of time $m+1$.

As a result, we divide our maximization problem by two problems:

1. Choose the optimal policy for the next moments of time for every initial distribution on $X_{m+1}$;

2. Choose the first step according to maximum reward and maximum value of the optimal policy assessment in the next time moments for initial distribution $\acute{\mu}$.

\begin{defin}[Derived model]
The model which is build of the model $Z^*$ by deletion $X_m$ and $A_{m+1}$ is called the \textbf{derived model} and it is denoted $\acute{Z^*}$.
\end{defin}

\begin{tver}[Fundamental equation]
$$\omega(x,\pi)=\sum\limits_{A(x)}\pi(a|x)\Big( q(a)+\acute{\omega}(p_a,\pi_a) \Big), \eqno (5.1) $$

where $p_a=p(\cdot|a), \pi_a(\cdot|\acute{h})=\pi(\cdot|ya\acute{h}),$

$a \in A_{m+1}, y=j(a), ~\acute{h}$ is a history in model
$\acute{Z^*}$.

The equation (5.1) is called \textbf{fundamental} and expresses the assessment $\omega$ of the random policy $\pi$ in model $Z^*$ in terms of the assessment $\acute{\omega}$ of some policies in the model $\acute{Z^*}$.
\end{tver}

\begin{flushleft}
\textbf{Proof}.
According to (4.2) we get
\end{flushleft}

$$\acute{\omega}(p_a,\pi_a)=\sum\limits_{X_{m+1}}p(y|a)\acute{\omega}(y,\pi_a) \eqno (5.2)$$

Let consider the spaces of ways $L$ and $\acute{L}$ in the models $Z^*$ and $\acute{Z^*}$.
Let $P^*$ is the probability distribution on $L$ according to the initial state $x$ and the policy $\pi$ and let $P^*_a$ is the probability distribution on $\acute{L}$ according to the initial distribution $p_a$ and the policy $\pi_a$.

\vspace{4mm}

According to (2.1) and (3.1) $\forall \acute{l}\in \acute{L}$ we get
$$I(xa\acute{l},x^*)=q(a)+I(\acute{l},x^*_{-1}) \eqno (5.3)$$
$$P^*(xa\acute{l},x^*)=\pi(a|x)P^*_a(\acute{l},x^*_{-1}) \eqno (5.4)$$
$$a\in A(x), x^*_{-1}=(x^*_{m+2},\ldots, x^*_{n}), (x^*_{m+1},x^*_{-1})=x^*.$$

Under the notations in (3.2) and (3.3) we get
$$\omega(x,\pi)=\sum\limits_{L}P^*(l,x^*)I(l,x^*) \eqno (5.5)$$
$$\acute{\omega}(p_a,\pi_a)=\sum\limits_{\acute{L}}P^*_a(\acute{l},x^*_{-1})I(\acute{l},x^*_{-1}) \eqno (5.6)$$

The measure $P^*(l,x^*)$ is nonzero only for ways which have the starting point $x$, i.e., for $xa\acute{l}$.
That is why by the substitution in (5.5) of the expression of $I(l,x^*)$ from (5.3) and the expression of $P^*(l,x^*)$ from (5.4), and according to (5.6) we get the fundamental equation (5.1). \textbf{\textit{Proposition 2 is proved.}}

\begin{remark}
The fundamental equation is correct even without \textbf{\textit{Assumption 1.}}
\end{remark}

\paragraph{6. Reducing the problem of the optimal decision to analogical problem for the derived model.}

From fundamental equation (5.1) it follows the following inequality

$$\omega(x,\pi)\leq\sup\limits_{A(x)}[q(a)+\acute{\omega}(p_a,\pi_a)]\leq\sup\limits_{A(x)}[q(a)+\acute{\nu}(p_a)] \eqno (6.1)$$

$\forall x\in X_m$ and for every $\pi$ ($\acute{\nu}$ which is the assessment of model $\acute{Z^*}$).

We denote $u(a)=q(a)+\acute{\nu}(p_a), (a\in A_{m+1})$ and call this value - \textbf{assessment of the action} $a$.

According to (4.3) and $\nu(x^*)=c(x^*)$ we get $u=U\acute{\nu}$ where operator $U$ transforms functions 
on the non-killed states on $X$ to the functions on $A$ and is given by

$$Uf(a)=q(a)+\sum\limits_{y}p(y|a)f(y)+\sum\limits_{y^*}p(y^*|a)c(y^*) \eqno (6.2)$$

where $y$ and $y^*$ are the non-killed states and the killed states, respectively.

Let the operator $V$ transforms the functions on $A$ into the functions on non-killed and non-terminal states on $X$ and satisfies
$$Vg(x)=\sup\limits_{a\in A(x)}g(a) \eqno (6.3)$$

Let us write the inequality (6.1) by using the operator $V$:
$$\omega(x,\pi)\leq Vu(x).$$

Then we consider $\sup\limits_{\pi}$ of the right and the left part of $\omega(x,\pi)\leq Vu(x)$ and we get 

$$\nu\leq Vu. \eqno(6.4)$$
\begin{remark}
Later we show the conditions which assure the equality in (6.4).
\end{remark}

\begin{defin}[Product of policies]
Let $\acute{\pi}$ be a killed policy in the model $\acute{Z^*}$ and to $x\in X_m$ is assigned some probability distribution $\gamma(\cdot|x)$ on $A_{m+1}$ which is concentrated on $A(x)$. When we choose on the first step an action $a$ and on all other steps we use the killed policy $\acute{\pi}$ then we get the killed policy $\pi$ in the model $Z^*$. This policy is called the \textbf{product of policies} $\gamma$ and $\acute{\pi}$ and is denoted by $\gamma\acute{\pi}$. It has the expression
\begin{displaymath}
\pi(\cdot|h)=\left \{
\begin{array}{ll}
\gamma(\cdot|x) &
\textrm{for} ~h=x\in X_m,\\
\acute{\pi}(\cdot|\acute{h}) & \textrm{for} ~h=xa\acute{h}.
\end{array}  \right .
\end{displaymath}
\end{defin}
\begin{tver}

Let $\pi=\gamma\acute{\pi}$ is a product of the killed policies $\gamma$ and $\acute{\pi}$. If $\acute{\pi}$ is uniform $\varepsilon'$-optimal for model $\acute{Z^*}$ then:

$$\nu=Vu. \eqno (6.4)$$

\end{tver}

\begin{flushleft}
\textbf{Proof.} The fundamental equation (5.1) for a product of policies has the following expression
\end{flushleft}

$$\omega(x,\gamma\acute{\pi})=\sum\limits_{A(x)}\gamma(a|x)\Big(q(a)+\acute{\omega}(p_a,\acute{\pi})\Big) \eqno (6.5)$$

Since $\acute{\pi}$ is $\varepsilon'$-optimal (it exists $\forall\ \varepsilon'\geq0$ according to \textit{Proposition 1.})
we get $\acute{\omega}(p_a,\acute{\pi})\geq\acute{\nu}(p_a)-\varepsilon'$, and according to appearance of $u$ equation (6.5) transforms to

$$\omega(x,\gamma\acute{\pi})\geq\sum\limits_{A(x)}\gamma(a|x)u(a)-\varepsilon'.$$

Lets consider the set

$$A_{\chi}(x)=\{a:a\in A(x), u(a)\geq Vu(x)-\chi \} ~(x\in X_m).$$

$A_{\chi}(x)$ is nonempty for all $\chi>0$. Let $\gamma(\cdot|x)$ be a probability distribution on $A(x)$ which is concentrated on $A_{\chi}(x)$.

Then
$$\sum\limits_{A(x)}\gamma(a|x)u(a)\geq Vu(x)-\chi.$$

Since $\varepsilon'+\chi\leq\varepsilon$ we get

$$\omega(x,\pi)\geq Vu(x)-\varepsilon, ~(x\in X_m). \eqno (6.6)$$

According to (6.4) and (6.6) \textit{\textbf{Proposition 3 is proved.}}

\vspace{4mm}

\begin{nasl}
The assessment $\nu$ of the model $Z^*$ is expressed in terms of the assessment $\acute{\nu}$ of the model $\acute{Z^*}$ in the following way: 

$$\nu=Vu, \ u=U\acute{\nu} \eqno (6.7)$$

where operators $U$ and $V$ are defined in (6.2) and (6.3);
\end{nasl}

\begin{nasl}
For all $\chi>0$ exists such $\psi(x):X_m\rightarrow A_{m+1}(x)$:

$$u(\psi(x))\geq\nu(x)-\chi \eqno (6.8)$$

Here $\gamma(\cdot|x)$ can be the distribution concentrated at one point $\psi(x)\in A_{\chi}(x)$.
\end{nasl}

\begin{nasl}
Let $\varepsilon'$ and $\chi$ be the arbitrary nonnegative numbers. If $\acute{\pi}$ is uniform $\varepsilon'$-optimal for the model $\acute{Z^*}$ and $\psi$ is such as in \textit{Corollary 3} then the killed policy $\psi\acute{\pi}$ is uniform $(\varepsilon'+\chi)$-optimal for the model $Z^*$.
\end{nasl}

\paragraph{7. Optimality equation. Method for the construction of simple optimal policies.}

Let assume that in our model $Z^* \ m=0$. Let consider the models $Z^*_0,Z^*_1,\ldots,Z^*_n$ where $Z^*=Z^*_0$ and $Z^*_t$ is a derived model of $Z^*_{t-1}$.
Let denote the assessments $\nu$ and $u$ of the model $Z^*_t$ as  $\nu_t$ and $u_{t+1}$($\nu_t$ on $X_t$, $u_{t+1}$ on $A_{t+1}$).
The reward function $q$ and the transition function $p$ we denote $q_t$ and $p_t$.

According to the results of \textit{section 6} we get

$$\nu_{t-1}=Vu_t, ~u_t=U\nu_t ~(1\leq t\leq n) \eqno (7.1)$$

where $$U_{t}f(a)=q_{t}(a)+\sum\limits_{y\in
X_t}p_t(y|a)f(y)+p_t(y^*|a)c(y^*), ~(a\in A_t, y^*\in X_t),$$

$$V_{t}g(x)=\sup\limits_{A(x)}g(a), ~(x\in X_{t-1}),$$

and $\nu_n=r.$

Equations (7.1) are called the \textbf{optimality equations}.
Let $T_t=V_{t}U_t$ then the optimality equations transform to

$$\nu_{t-1}=T_{t}\nu_t. \eqno (7.\acute{1})$$

From (7.1),(7.$\acute{1}$) and the condition $\nu_n=r$ we calculate $\nu_n,\nu_{n-1},\ldots,\nu_0$.
Then we choose the action $\psi_t(x):X_{t-1}\rightarrow A_{t}(x)$ for which holds

$$u_t(\psi_t)\geq\nu_{t-1}-\chi_t.\emph{\emph{\emph{\emph{}}}} \eqno (7.2)$$

$\forall t=1,2,\ldots,n$ and for all nonnegative $\chi_1,\chi_2,\ldots\chi_n$.

According to \textit{Corollary 3} of \textit{Proposition 3} the simple policy $\varphi=\psi_{1}\psi_{2}\ldots\psi_n$ is uniform $\varepsilon$-optimal for the model $Z^*=Z^*_0$ and $\varepsilon=\sum\limits_{i=1}^{n}\chi_i$.
The equation (7.2) can be rewritten as

$$T_{\psi_t}\nu_t\geq\nu_{t-1}-\chi_t, \eqno (7.\acute{2})$$

where the operator $T_{\psi_t}$ transforms functions on $X_t$ to functions on $X_{t-1}$ in the following way

$$T_{\psi_t}f(x)=q_{t}[\psi_{t}(x)]+\sum\limits_{X_t}p(y|\psi_{t}(x))f(y)+p_t(y^*|a)c(y^*). \eqno (7.3)$$

\begin{tver}
Let $\pi$ be an arbitrary killed policy in the derived model $Z^*_k \
(k=1,2,\ldots,n)$ and let $\psi_t:X_{t-1}\rightarrow A_t(x) \ (t=1,2,\ldots,k)$ are arbitrary too then 

$$\omega_0(x,\psi_{1}\psi_{2}\ldots\psi_{k}\pi)=T_{\psi_{1}}T_{\psi_{2}}\ldots T_{\psi_k}\omega_{k}(x,\pi), \eqno (7.4)$$

\end{tver}

\begin{flushleft}
\textbf{Proof}. 
It follows from the fundamental equation (5.1), formulas (5.2), (7.3) and the mathematical induction.
\end{flushleft}
\begin{remark}
It follows from (7.4): the result will not change if our decision process is killed at the moment of time $k$ and the terminal reward as the assessment of policy $\pi$ is taken.
\end{remark}
\begin{remark}
If we can choose $\psi_t$ with $\chi_t=0$ in $(7.\acute{2})$ $\forall t=1..n$ then the simple policy $\varphi=\psi_1\ldots\psi_n$ is called uniform optimal.
\end{remark}

\paragraph{8. The sufficiency of the simple policies for countable models.}

The question arises: do we lose something by using only simple policies?
The previous result can not give us the answer. It only makes our losses indefinitely small.

\begin{teorema}[Sufficiency of the simple policies]
Let $\mu$ is a fixed initial distribution and let $\pi$ is a arbitrary killed policy then there exists $\varphi$-simple policy such that

$$\omega(\mu,\pi)\leq\omega(\mu,\varphi).\eqno(8.1)$$

\end{teorema}
\begin{flushleft}
\textbf{Proof.} It follows from \textit{Proposition 5} and \textit{Proposition 6}.
\end{flushleft}

\begin{tver}
For all $\mu$ and for all killed policies $\pi$ there exists the Markov policy $\theta$ such that

$$\omega(\mu,\theta)=\omega(\mu,\pi)\eqno(8.2)$$
These two policies are called \textbf{equivalent}.
\end{tver}
\begin{tver}
For all Markov policies $\theta$ there exists a simple policy $\varphi$ such that
$$\omega(\mu,\varphi)\geq\omega(\mu,\theta)\eqno(8.3)$$
We say that $\varphi$ \textbf{dominates} $\theta$ \textbf{uniformly}.
\end{tver}
\begin{flushleft}
\textbf{Proof.}(\textit{Proposition 5}). Let $\theta$ is Markov policy and
\end{flushleft}

$$\theta(a|x)=\mathbb{P^*}\{a_t=a|x_{t-1}=x\}=\frac{\mathbb{P^*}\{x_{t-1}a_t=xa\}}{\mathbb{P^*}\{x_{t-1}=x\}}\eqno(8.4)$$
\begin{center}
$(a\in A_t,~~x\in X_{t-1},~~m+1\leq t\leq n)$,
\end{center}
where $\mathbb{P^*}$ is a probability measure in the space of ways $L$ which is assigned  to the initial distribution $\mu$ and to the policy $\pi$.
\begin{remark}
The expression on the right side of (8.4) makes no sense for $\mathbb{P^*}\{x_{t-1}=x\}=0$.
So, for such $x$(in particular for killed states) we choose the arbitrary distribution on $A(x)$ instead of $\theta(\cdot|x)$.
\end{remark}
Let $\mathbb{Q^*}$ denotes a probability distribution on space $L$ which is assigned  to the initial distribution $\mu$ and to the killed Markov policy $\theta$.

The distribution $\mathbb{Q^*}$ does not match with $\mathbb{P^*}$ in the general case, but it is enough for proving (8.2) if any of $x_m,a_{m+1},\ldots,a_n,x_n$ and $x^*_{m+1},x^*_{m+2},\ldots,x^*_{n}$ has the same probability distribution according to measures $\mathbb{P^*}$ and $\mathbb{Q^*}$.

The following assertion holds

$$\omega(\mu,\pi)=\sum\limits_{t=m+1}^n\mathbb{P^*}q(a_t)+\sum\limits_{t=m+1}^n\mathbb{P^*}c(x^*_t)+\mathbb{P^*}r(x_n),$$
$$\omega(\mu,\theta)=\sum\limits_{t=m+1}^n\mathbb{Q^*}q(a_t)+\sum\limits_{t=m+1}^n\mathbb{Q^*}c(x^*_t)+\mathbb{Q^*}r(x_n).$$

We shall use the mathematical induction to prove this.

The \textbf{basis} of induction: (8.2) holds for $x_m$ because $\mathbb{P^*}=\mathbb{Q^*}=\mu$.

The \textbf{induction hypothesis}: let (8.2) holds for $x_{t-1}$. Let's check it for $a_t$.

Since $\theta$ is a killed Markov policy then 

$$\mathbb{Q^*}\{x_{t-1}a_t=xa\}=\mathbb{Q^*}\{x_{t-1}=x\}\theta(a|x),~~(a\in A_t,~x\in X_{t-1}).\eqno(8.5)$$

Hence, from (8.4) and (8.5) we get

$$\mathbb{P^*}\{a_t=a\}=\sum\limits_{x\in X_{t-1}}\mathbb{P^*}\{x_{t-1}a_t=xa\}=\sum\limits_{x\in X_{t-1}}\mathbb{P^*}\{x_{t-1}=x\}\theta(a|x)=$$
$$=\sum\limits_{x\in X_{t-1}}\mathbb{Q^*}\{x_{t-1}=x\}\theta(a|x)=\sum\limits_{x\in X_{t-1}}\mathbb{Q^*}\{x_{t-1}a_t=xa\}=\mathbb{Q^*}\{a_t=a\}.$$

So, our proposition holds for $a_t$.

The \textbf{induction hypothesis}: let (8.2) holds for $a_{t}$. Let show it for $x_t$. 

From the definition of the transition function we get

$$\mathbb{P^*}\{a_{t}x_t=ax\}=\mathbb{P^*}\{a_t=a\}p(x|a),\eqno(8.6)$$

$$\mathbb{Q^*}\{a_{t}x_t=ax\}=\mathbb{Q^*}\{a_t=a\}p(x|a).\eqno(8.7)$$

From (8.6) and(8.7) it follows

$$\mathbb{P^*}\{x_t=x\}=\sum\limits_{a\in A_t}\mathbb{P^*}\{a_{t}x_t=ax\}=\sum\limits_{a\in A_t}\mathbb{P^*}\{a_t=a\}p(x|a)=$$
$$=\sum\limits_{a\in A_t}\mathbb{Q^*}\{a_t=a\}p(x|a)=\sum\limits_{a\in A_t}\mathbb{Q^*}\{a_{t}x_t=ax\}=\mathbb{Q^*}\{x_t=x\},~~(x\in X_t).$$

\textbf{\textit{Proposition 5 is proved.}}

\vspace{4mm}

\begin{flushleft}
\textbf{Proof.}(\textit{Proposition 6.})
For proving this proposition we need the following lemma.
\end{flushleft}

\begin{lema}
Let $f$ is a arbitrary function and let $\nu$ is a arbitrary probability distribution on countable space $E$.

If $\nu f<+\infty$ then the set $\Gamma=\{x:f(x)\geq\nu f\}$ has a positive measure $\nu$, namely
$$\nu(\Gamma)>0$$
(See proof in [1]).
\end{lema}

According to (4.2) the condition (8.3) is equal to

$$\omega(x,\varphi)\geq\omega(x,\theta),~~\forall x\in X_m.$$

Let separate the killed Markov policy $\theta$ by a product of the policies $\theta=\gamma\theta'$ where $\gamma$ is the restriction of $\theta$ on $X_m$ and $\theta'$ is the restriction of $\theta$ on $X_{m+1}\bigcup X_{m+2}\ldots \bigcup X_n$.

According to the fundamental equation (5.1) it holds

$$\omega(x,\theta)=\gamma_x f,$$

where $\gamma_x(\cdot)=\gamma(\cdot|x)$ is the probability distribution on $A(x)$, 

and $f(a)=q(a)+\omega'(p_a,\theta'),~~(a\in A_{m+1}).$

Since \textbf{\textit{Lemma 1}} for $\tilde{A}(x)\subset A(x)$ it follows $\gamma_x(\tilde{A}(x))>0$,
where $\tilde{A}(x)=\{a:f(a)\geq\gamma_x f=\omega(x,\theta)\}$. As a result, $\tilde{A}(x)$ is nonempty.
If $\psi(x)$ is an arbitrary point of $\tilde{A}(x)$ then $f(\psi(x))\geq\omega(x,\theta)$.
But since the fundamental equation (5.1) we get $f(\psi(x))=\omega(x,\psi\theta')$ and

$$\omega(x,\psi\theta')\geq\omega(x,\theta).$$

Let assume that condition (8.3) holds for the derived model $\acute{Z^*}$. Then exists a simple policy $\varphi'$ in $\acute{Z^*}$ which uniformly dominates the killed Markov policy $\theta'$. According to the fundamental equation (5.1) and our assumption we get

$$\omega(x,\psi\varphi')=q(\psi(x))+\omega'(p_{\psi(x)},\varphi')\geq q(\psi(x))+\omega'(p_{\psi(x)},\theta')=\omega(x,\psi\theta')\geq\omega(x,\theta).$$

In the model $Z^*$ simple policy $\varphi=\psi\varphi'$ dominates $\theta$ uniformly. 
Finally, (8.3) holds for model $Z^*$ too.

\textbf{\textit{Proposition 6. is proved.}}

\paragraph{9. Markovian property.}

Let $0<k<n$, let use the killed policy $\rho$ on the interval $[0,k]$ and killed policy $\pi$ on the interval $[k,n]$.
Doing analogically to \textbf{\textit{Definition 15}} we can say that policy $\rho\pi$ is used.

\begin{tver}
Let $L_0$ is the space of ways on the interval $[0,n]$, let $L_k$ is the space of ways on the interval $[k,n]$ and let $P^{*\rho\pi}_x$ is the probability distribution which is assigned to the initial state $x$ and to the killed policy $\rho\pi$, and analogically $P^{*\pi}_y$ is the probability distribution on $L_k$.

Then $\forall \xi=\xi(x_{k}a_{k+1}\ldots
x_n)$ on $L_k$ holds
 
$$E^{*\rho\pi}_x\xi=E^{*\rho}_{x}[E^{*\pi}_{x_{k}}\xi]. \eqno (9.1)$$

\end{tver}

\begin{flushleft}
\textbf{Proof}. $\forall l=y_{0}b_{1}\ldots
b_{k}y_{k}b_{k+1}\ldots y_n$ according to (3.1)
\end{flushleft}

$$P^{*\rho\pi}_x(y_{0}b_{1}\ldots y_n)=P^{*\rho}(cy_k)P^{*\pi}_{y_k}(y_{k}d), \eqno (9.2)$$

where $c=y_{0}b_{1}\ldots b_{k}$, $d=b_{k+1}\ldots y_n$. 
Any function $\xi$ on the space $L_k$ can be interpreted on $L_0$ like function which does not depend on $x_{0}a_{1},\ldots,a_{k}$.
That is why we multiply the both sides of (9.2) by $\xi(y_{k}d)$ and sum up over all ways

$$E^{*\rho\pi}_{x}\xi=\sum\limits_{cy_k}P^{*\rho}_{x}(cy_k)\sum\limits_{d}P^{*\pi}_{y_k}(y_{k}d)\xi(y_{k}d). \eqno (9.3)$$

But $P^{*\pi}_{y_k}(yd)=0 ~\textrm{for} ~y\neq y_k$ and it follows

$$\sum\limits_{d}P^{*\pi}_{y_k}(y_{k}d)\xi(y_{k}d)=\sum\limits_{yd}P^{*\pi}_{y_k}(yd)\xi(yd)=F(y_k). \eqno (9.4)$$

By substitution in (9.3) the expression from (9.4) and according to
$\sum\limits_{cy_k}P^{*\rho}_{x}(cy_k)F(y_k)=E^{*\rho}_{x}F(x_k)$,
we get (9.1). \textit{\textbf{Proposition 7 is proved.}}

\begin{flushleft}
\textbf{Corollary 1.(Markovian property)}
Let $\nu(y)=P^{*\rho}_{\mu}\{x_{k}=y\} \ (y\in X_k)$ then 
$\forall \mu$
\end{flushleft}
$$E^{*\rho\pi}_\mu\xi=E^{*\rho}_{\mu}[E^{*\pi}_{x_{k}}\xi].$$

In particular $$E^{*\rho\pi}_{\mu}\xi(x_{k}a_{k+1}\ldots
x_n)=E^{*\pi}_{\nu}\xi(x_{k}a_{k+1}\ldots x_n), \eqno (9.5)$$

It follows form (9.1) and $\sum\limits_{y\in
X_k}\nu(y)P^{*\pi}_{y}\xi=E^{*\pi}_{\nu}\xi$.

The formula (9.5) shows that the probability distribution for a part of the trajectory does not depend on the distribution $\mu$ and policy $\rho$ on the interval  $[k,n]$. Namely, the probability forecast of the "future"$(\xi)$ depends not on the "past" $(\mu,\rho)$, but only on the "present" $(\nu)$.
Actually, it is already the \textbf{Markovain property}.

Let use the Markovian property for the assessment of a killed policy $\rho\pi$ on the intervals $[0,k]$ and $[k,n]$.
Instead of $\xi$ we take $\xi=\sum\limits_{t=k+1}^{n}[q(a_t)+c(x^*_t)]+r(x_n)$ and by substituting in (9.5) we get

$$\omega(\mu,\rho\pi)=\sum\limits_{t=1}^{k}E^{*\rho\pi}_{\mu}[q(a_t)+c(x^*_t)]+\omega(\nu,\pi)=
\sum\limits_{t=1}^{k}E^{*\rho}_{\mu}[q(a_t)+c(x^*_t)]+\omega(\nu,\pi). \eqno (9.6)$$

The summation in (8.6) expresses the assessment $\omega(\mu,\rho)$ of policy $\rho$ for a zero terminal reward, namely, $\omega(\mu,\rho\pi)=\omega(\mu,\rho)+\omega(\nu,\pi)$.

There is also another interpretation of (9.6). According to (4.2) and $\nu(y)=P^{*\rho}_{\mu}\{x_{k}=y\} \ (y\in X_k)$ we get
$$\omega(\nu,\pi)=\sum\limits_{y}\nu(y)\omega(y,\pi)=E^{*\rho}_{\mu}\omega(x_k,\pi),$$
$$\omega(\mu,\rho\pi)=E^{*\rho}_{\mu}[\sum\limits_{t=1}^{k}q(a_t)+\omega(x_k,\pi)]. \eqno (9.7)$$

Hence, the assessment of killed policy $\rho\pi$ is equal to the assessment of the killed policy $\rho$ with the terminal reward $\omega(\cdot,\pi)$ at the moment of time $k$.

\paragraph{10. Dynamic programming principle.}

Let $Z^*$ be the model on the interval $[0,n]$ and let $0\leq s<t\leq n$.
Let $Z^*_{s,t}[f]$ denotes the model which is taken from the model $Z^*$ by restriction of the interval $[0,n]$ to $[s,t]$. We define the terminal reward $f$ at the moment of time $t$. Moreover, denote $\nu^t_s[f]$ as the assessment of the model $Z^{*t}_{s}$ with the terminal reward $f$. Obviously,  $\nu^t_s[f]=(VU)^{t-s}f=T^{t-s}f$ on $X$.

Since $\forall t\in [0,n]$ it holds

$$\nu^n_0[r]=\nu^t_0[\nu^n_t[r]] ~\textrm{on} ~X_0 ~(r ~\textrm{on} ~X_n). \eqno (10.1)$$

The equation (10.1) is equivalent to the optimality equations (7.1) and the condition $\nu^n=r$.
It is called the \textit{\textbf{Dynamic programming principle}} and it means that for the optimization of the decision on the interval $[0,n]$ with terminal reward $r$ we must first optimize the decision on interval $[t,n]$(with such terminal reward) and then optimize the decision on the interval $[0,t]$ with terminal reward $\nu^n_t[r]$.

In particular according to (9.1) it follows that if $\pi''$ is a uniform $\varepsilon$-optimal killed policy for $Z^{*n}_t$ with terminal reward $r$ and 
$\pi'$ is a uniform $\varepsilon$-optimal policy for $Z^{*t}_0$ with the terminal reward $\nu^n_t[r]$ then the
killed policy $\pi=\pi''\pi'$ has the assessment $\nu^n_0[r]$ and is uniform $\varepsilon$-optimal for the model $Z^{*n}_0$(with terminal reward $r$).

\begin{center}
\textbf{References}
\end{center}

\begin{flushleft}
[1]. E.B. Dynkin, A.A. Yushkevich, \textit{Markov Decision Processes}, M., (1975), 334 p. \linebreak(Russian)

[2]. E.A. Feinberg, A. Shwartz,  \textit{Introduction}, Handbook of Markov Decision Processes, Kluwer, (2002)(565 pages), pp.1-17.(English)

[3]. A.G. Pakes, \textit{Killing and Resurrection of Markov Processes}, Stochastic Models, V.13, I.2, (1997), pp.255-269.(English)

[4]. R.E. Bellman, \textit{Dynamic Programming}, Izdatelstvo inostrannoj literatury, (1960), 400 p.(Russian)
\end{flushleft}

\begin{flushleft}

\textbf{Nestor R. Parolya, Yaroslav I. Yeleyko}

Ivan Franko National University of Lviv

1, Universytetska str., 79000, Lviv, Ukraine

Tel.: \,\,(8032) 239 45 31 \,\,(off.)
\end{flushleft}

\end{document}